\documentclass[onefignum,onetabnum]{siamart250211}

\usepackage{lipsum}
\usepackage{amsfonts,amssymb,mathtools}
\usepackage{hyperref}
\usepackage{graphicx}
\usepackage{placeins}
\usepackage{epstopdf}
\usepackage{multirow}
\usepackage{comment}
\usepackage{algorithm}
\usepackage{algpseudocode}
\usepackage{enumitem}
\usepackage{booktabs}
\usepackage{color}

\ifpdf
  \DeclareGraphicsExtensions{.eps,.pdf,.png,.jpg}
\else
  \DeclareGraphicsExtensions{.eps}
\fi

\definecolor{red2}{RGB}{204,0,0}
\definecolor{blue2}{RGB}{0,103,165}
\definecolor{darkviolet}{RGB}{148,0,211}
\hypersetup{colorlinks,linkcolor=[RGB]{0,103,165},citecolor=[RGB]{180,0,0}}
\usepackage[nocompress]{cite}

\newsiamremark{remark}{Remark}
\newsiamremark{hypothesis}{Hypothesis}
\crefname{hypothesis}{Hypothesis}{Hypotheses}
\newsiamthm{claim}{Claim}

\def\be{\begin{equation}}
\def\ee{\end{equation}}

\def\Pi{\mathbf{\Phi}}

\def\Ut{\widetilde{U}}
\def\yt{\widetilde{y}}
\def\tPhi{\widetilde{\Phi}}

\newcommand{\argmin}{\operatornamewithlimits{argmin}}

\newcommand{\refm}[1]{\textcolor{darkviolet}{#1}}

\newcommand{\spanop}{\operatorname{span}}
\newcommand{\dist}{\operatorname{dist}}
\newtheorem{assumption}{Assumption}

\usepackage{amsopn}
\DeclareMathOperator{\diag}{diag}

\headers{Delay Closure for the Heat Equation}{Xiaoting Han, Dongbin Xiu, and Ludmil Zikatanov}

\title{Delay Closure of Observable Dynamics for the Heat Equation}
\hypersetup{
  pdftitle={Delay Closure of Observable Dynamics for the Heat Equation},
  pdfauthor={Xiaoting Han, Dongbin Xiu, and Ludmil T. Zikatanov}
}

\author{Xiaoting Han\thanks{Department of Mathematics, The Ohio State University, Columbus, OH 43210, USA. E-mail: \texttt{han.2051@osu.edu}, \texttt{xiu.16@osu.edu}. Funding: This work was partially supported by AFOSR FA9550-24-1-0237.} \and Dongbin Xiu\footnotemark[1]\and Ludmil T. Zikatanov\thanks{Department of Mathematics, Penn State, University Park, PA 16802, USA. E-mail: \texttt{ludmil@psu.edu}. Funding: The work of Zikatanov was supported in part by the Center for Computational Mathematics and its Applications at Penn State and by the Center for Advanced Studies at the Norwegian Academy of Science and Letters, Norway under the program Mathematical Challenges in Brain Mechanics (2025/2026).} }

\begin{document}

\maketitle

\begin{abstract}
We study finite-delay representations of point observables, which are solution values at selected spatial locations, for spatially discretized heat equations. We show that the observables satisfy closed dynamical systems whose evolution depends only on their own history, without requiring solution values elsewhere in the domain. Using semigroup properties, we analyze low-rank approximate closures that can use substantially fewer delay terms than the exact representation. We consider two computational approaches: direct least-squares fitting of delay coefficients from observable trajectories and shift-invert Krylov approximation of the semigroup using preconditioned linear solves. Numerical experiments for constant- and variable-coefficient heat equations in one, two, and three dimensions illustrate the effects of memory length and sensor information on prediction. The results are consistent with the theoretical analysis.
\end{abstract}



\section{Introduction}
In many prediction problems, the quantity of interest is an observable of a dynamical system rather than the complete state. For a diffusion process, this observable may be the solution at one or several sensor locations. Although the full state satisfies a closed governing equation, a point observation generally does not form a closed system, because its evolution depends on state variables at other locations.

The focus of this paper is on closed representations of observable dynamics, that is, dynamical laws for observables without explicit dependence on the full state. More specifically, we study such representations for the heat equation after spatial discretization. Our starting point is the linear system
\[
    \dot U=-A_hU,\qquad y=b^TU,
\]
where $A_h\in\mathbb{R}^{N\times N}$ is a symmetric positive-definite discrete diffusion matrix, $U\in\mathbb{R}^N$ is the state, and $b\in\mathbb{R}^N$ defines a linear observation. Liu and Xiu~\cite{LiuXiu_JMCMC26} establish the existence of an exact finite-delay representation for the linear observable $y$. Although useful mathematically, this representation may require as many as $N$ delays, where $N$ is determined by the spatial discretization and can be large. On the other hand, delay representations have been used extensively to model observable dynamics through flow map learning (FML)~\cite{Churchill_2023,FuChangXiu_JMLMC20,wu2020data,chen2022deep}. For a variety of linear and nonlinear problems, the number of delays has been modest and has shown little, if any, dependence on the state dimension $N$. The exact algebraic closure result of \cite{LiuXiu_JMCMC26} therefore does not explain why a small number of delays can produce accurate observable predictions on a fine spatial mesh.

In this paper, we consider the semigroup over a fixed time lag $\tau>0$,
\[
    \Phi_\tau=e^{-\tau A_h}.
\]
Using an approach similar to that of~\cite{LiuXiu_JMCMC26}, we establish an exact delay representation based on the semigroup Krylov space generated by $b,\Phi_\tau b,\ldots$. Although this exact representation may still require $N$ delays, the semigroup formulation allows us to exploit the damping of high-frequency modes over time. For the heat equation, this damping supports accurate approximations of the observable dynamics with substantially fewer delays.

The analysis also distinguishes smoothing time from memory length. A shifted Krylov sequence begins after a smoothing period, often called burn-in, which reduces the contribution of high-frequency components. Discarding an initial transient is a common procedure in FML~\cite{Churchill_2023}. The present analysis quantifies the role of this smoothing for the heat equation.

For the numerical construction of a delay closure, we consider two approaches. The first constructs a low-dimensional representation of semigroup actions using shift-invert Krylov methods. Spectral-transformation methods~\cite{EricssonRuhe1980}, Krylov approximation of the matrix exponential~\cite{Saad1992,HochbruckLubich1997,EshofHochbruck2006}, and multilevel preconditioners~\cite{BramblePasciakXu1990,BrandtMcCormickRuge1982,RugeStuben1987} provide the numerical ingredients. The second is the data-driven approach used in existing FML studies, where linear delay coefficients are fitted by least squares to short observable segments sampled from trajectories with independent initial conditions.

The numerical examples first examine semigroup approximation and effective rank under mesh refinement. They then assess autonomous data-driven prediction for constant and spatially varying diffusivities in one, two, and three dimensions. Further sensor studies examine the number and placement of observations in one dimension and point-observable prediction at interior and Neumann-boundary locations in two and three dimensions.

The paper is organized as follows. Section~\ref{sec:setup} specifies the semidiscrete setting. Section~\ref{sec:delay} derives exact sampled-time closure using semigroup Krylov spaces, and Section~\ref{sec:low-rank} develops low-rank approximate closure via spectral-tail and smoothed Krylov approximation estimates. Section~\ref{sec:computation} describes operator-based and data-driven construction of closure models. Section~\ref{sec:examples} presents the numerical results.

\section{Problem Setup}\label{sec:setup}

Consider the heat equation
\begin{equation}\label{heat}
  u_t = \nabla\cdot\left(a(x)\nabla u\right), \qquad
  0<a_{\min}\leq a(x)\leq a_{\max}<\infty,
\end{equation}
on the unit hypercube $\Omega=(0,1)^d$, $d\in\{1,2,3\}$, with
homogeneous boundary conditions and a nonempty Dirichlet portion of
the boundary.
We define the observable as the solution at a fixed location,
\begin{equation}\label{obs}
  y(t)=u(x^*,t),\qquad x^*\in\overline\Omega,
\end{equation}
where $x^*$ lies in the interior or on a Neumann boundary. Points on a homogeneous Dirichlet boundary are excluded because their observable dynamics are trivial. Our goal is to understand and compute the dynamics of this observable without using the solution at other locations during prediction.

Specifically, we focus on the semidiscrete system obtained by spatial discretization. Let $N>1$ be the number of independent discrete unknowns and $h>0$ the characteristic mesh size. We consider
\begin{equation}\label{Ut}
  \dot U=-A_hU,
\end{equation}
where $U\in\mathbb R^N$ is the discrete solution vector and $A_h\in\mathbb R^{N\times N}$ is a symmetric positive-definite matrix. The discrete approximation to the observable~\eqref{obs} is represented by
\begin{equation}\label{y}
  y=b^TU,\qquad b\in\mathbb R^N\setminus\{0\}.
\end{equation}
If $x^*$ is a grid point, say $x^*=x_i$, then $b=e_i$, the $i$th canonical basis vector. Otherwise, the nonzero entries of $b$ are interpolation weights associated with nearby grid points and define a local approximation to $u(x^*,t)$.

All vector norms and induced matrix norms are Euclidean unless specified otherwise. The symmetry and positive definiteness of $A_h$ imply
\begin{equation}\label{eq:state-contraction}
  \left\|e^{-tA_h}\right\|_2=e^{-t\lambda_{h,1}}\leq1,\qquad t\geq0,
\end{equation}
where $\lambda_{h,1}>0$ is its smallest eigenvalue.

\section{Exact Delay Closure for Dynamical Systems} \label{sec:delay}
In this section, we discuss exact delay representations of observable dynamics
for a general linear dynamical system, not restricted to the semidiscrete heat equation~\eqref{Ut}. Consider
\begin{equation} \label{Dyn}
  \dot{U} = -A U, \qquad y=b^T U,
\end{equation}
where $A\in\mathbb R^{N\times N}$ needs not be symmetric or positive definite and $b\in\mathbb R^N\setminus\{0\}$. We first discuss the exact delay representation from~\cite{LiuXiu_JMCMC26}, termed algebraic delay closure, and then establish a new semigroup delay closure.

\subsection{Algebraic Delay Closure} \label{sec:algebraic}

For the system~\eqref{Dyn}, consider the Krylov space
\begin{equation} \label{Kry0}
  \mathcal{K}_m(A^T,b) \coloneqq \operatorname{span}\left\{b,A^T b,\left(A^T\right)^2b,\dots, \left(A^T\right)^{m-1}b\right\}, \quad m\geq 1.
\end{equation}
Let the limiting dimension of the Krylov space be $n$, i.e.,
\begin{equation}
  n \coloneqq \lim_{m\to\infty} \dim (\mathcal{K}_m(A^T,b)) \leq N.
\end{equation}
Theorem~1 of~\cite{LiuXiu_JMCMC26} gives an exact finite-delay representation for $\dot y$ for a suitable choice of distinct delays,
\begin{equation} \label{ydot0}
\dot{y}(t) = \sum_{k=1}^{n} \alpha_k y(t-\tau_k), \qquad
0\leq \tau_1<\cdots<\tau_n,
\end{equation}
where $\alpha_1,\dots,\alpha_n$ are the unique coefficients for this choice of delays. This delay representation is determined by the observable algebraic structure of the pair $(A,b)$ and the chosen delays and will be referred to as ``algebraic closure'' hereafter.

While the algebraic closure is exact, its practical utility can be
limited. Consider the one-dimensional heat
equation~\eqref{Ut}--\eqref{y} with constant diffusivity
$a(x)\equiv1$, discretized by second-order central differences on a
uniform mesh of size $h$. The discrete diffusion matrix $A_h$ has the
interior stencil $h^{-2}\,\mathrm{tridiag}(-1,2,-1)$, with boundary
rows modified as required by the boundary conditions. For a point
observation $b=e_{i^*}$, successive applications of $A_h$ spread the
support of the Krylov vectors to neighboring nodes. When the
observable Krylov space has full dimension, $n=N$, the
representation~\eqref{ydot0} uses $N$ delays.
Since accurate spatial discretization can require $N\gg1$, the exact
algebraic closure can be an expensive dynamical model. It also does
not explain the modest number of delays observed in FML models, which have been
shown to be accurate for a variety of dynamical systems \cite{Churchill_2023,FuChangXiu_JMLMC20,wu2020data,chen2022deep}.

\subsection{Semigroup Delay Closure}
\label{sec:semigroup}

We now establish an alternative delay representation for~\eqref{Dyn} using a Krylov space generated by its semigroup.
Consider discrete times $t_k=k\tau$, $k=0,1,\ldots$, with a constant step $\tau>0$, and let $U_k\coloneqq U(t_k)$, where $U(t)$ solves~\eqref{Dyn}. The 
time-$\tau$ flow map $\Phi_\tau$ of \eqref{Dyn} is
\begin{equation}\label{Phi}
  \Phi_\tau=e^{-\tau A}:\mathbb{R}^N\to\mathbb{R}^N,
\end{equation}
and
\begin{equation} \label{Uk_yk}
U_{k+1} = \Phi_\tau U_k, \qquad y_k = b^TU_k.
\end{equation}
Applying $\Phi_\tau$ repeatedly, we obtain 
\begin{equation} \label{y-map}
    y_{k+\ell} = b^T \Phi_\tau^\ell U_k = \left(\left(\Phi_\tau^T\right)^\ell b\right)^T U_k, \qquad \ell \in \mathbb{Z}_{\ge 0}.
\end{equation}
With this notation in hand, we define the semigroup Krylov space as follows.
\begin{definition}[Semigroup Krylov space] \label{def:Kry}
  For the dynamical system~\eqref{Dyn} with the time-$\tau$ flow map~\eqref{Phi}, the semigroup Krylov space of order
  $m\geq 1$ is defined as
  \begin{equation} \label{Krylov}
     \mathcal K_m(\Phi^T_\tau,b)=\operatorname{span}\left\{b,\Phi_\tau^T b,\dots,\left(\Phi_\tau^T\right)^{m-1}b\right\}, \qquad m=1,2,\dots.
    \end{equation}
  \end{definition}
  The Krylov dimensions are nondecreasing and satisfy
  \[
    1\leq\dim(\mathcal K_m(\Phi_\tau^T,b))\leq\min\{m,N\},\qquad m\geq1.
  \]
  This yields the following result.
  \begin{theorem}[Semigroup closure] \label{thm:semigroup-closure}
Let $n\geq 1$ be the limiting dimension of the semigroup Krylov space \eqref{Krylov},
\begin{equation} \label{n_rank}
  n \coloneqq \lim_{m\to\infty} \dim (\mathcal{K}_m(\Phi^T_\tau,b)) \leq N.
\end{equation}
Then the observable of~\eqref{Dyn}, sampled over the discrete times \eqref{Uk_yk}, satisfies
 \begin{equation}\label{yk-closure}
    y_{k+1} = \sum_{\ell=0}^{n-1} {c}_\ell y_{k-\ell}, \qquad k\geq n-1,
\end{equation}
where $c_0, \dots, c_{n-1}$ are constants that are not all zero. 
\end{theorem}
\begin{proof}
  The Krylov dimensions are nondecreasing, and the limiting dimension $n$ is first attained at $\mathcal K_n(\Phi_\tau^T,b)$. Consequently, there exist constants $d_0,\ldots,d_{n-1}$, not all zero, such that
   \begin{equation}  \label{Phi-span}
     \left(\Phi_\tau^T\right)^n b = \sum_{\ell=0}^{n-1} d_\ell \left(\Phi_\tau^T\right)^\ell b.
   \end{equation}
   Taking the transpose and multiplying by $U_k$, we obtain
   \begin{equation*} \label{bPhiU}
b^T\Phi_\tau^n U_k = \sum_{\ell=0}^{n-1} d_\ell\, b^T \Phi_\tau^\ell U_k.
\end{equation*}
Using~\eqref{y-map}, we obtain
  \begin{equation} \label{yk1-closure}
y_{k+n} = \sum_{\ell=0}^{n-1} d_\ell y_{k+\ell}.
\end{equation}
Reindexing yields~\eqref{yk-closure} with
\begin{equation}\label{eq:coefficient-reversal}
    c_\ell=d_{n-1-\ell},\qquad \ell=0,\ldots,n-1.
\end{equation}
\end{proof}

The semigroup closure~\eqref{yk-closure} is exact in discrete time. As with the algebraic closure~\eqref{ydot0}, it can require $n=N$ delay terms for exactness.

\section{Low-Rank Semigroup Closure} \label{sec:low-rank}

For the remainder of the paper, we focus on the semidiscrete heat equation \eqref{Ut}, where the matrix $A_h$ is symmetric positive definite. Consequently, its time-$\tau$ flow map $\Phi_\tau=e^{-\tau A_h}$ is also symmetric. We therefore omit the transpose on $\Phi_\tau$ in what follows.

Following \eqref{Phi-span} in the proof of Theorem~\ref{thm:semigroup-closure}, suppose that, for some integer $1\leq r<n$,
\begin{equation} \label{Phi-approx}
  \left\|
  \Phi_\tau^r b
  -
  \sum_{\ell=0}^{r-1} d_\ell \Phi_\tau^\ell b
  \right\|_2
  \le \varepsilon,
\end{equation}
then the argument leading to~\eqref{yk1-closure} gives an order-$r$ approximate closure satisfying
\begin{equation} \label{apprx-closure}
  \left|
  y_{k+r}
  -
  \sum_{\ell=0}^{r-1} d_\ell y_{k+\ell}
  \right|
  \le
  \varepsilon \|U_k\|_2 \leq \varepsilon \|U_0\|_2,
\end{equation}
where the final inequality follows from contractivity~\eqref{eq:state-contraction}. 
Thus, the approximate
semigroup Krylov relation \eqref{Phi-approx} directly yields an approximate discrete-time delay closure with the controllable error bound \eqref{apprx-closure}.

\subsection{Visibility of Spectral Components}

The properties of the low-rank semigroup closure \eqref{apprx-closure} depend critically on the spectrum of the discrete diffusion matrix $A_h$ in~\eqref{Ut}.
Let $\lambda_{h,1}\leq\cdots\leq\lambda_{h,N}$ be the eigenvalues of $A_h$, counted with multiplicity. For a family of uniform-mesh discretizations, we assume that there is a constant $\widetilde{c}>0$, independent of $h$ and $j$, such that
\begin{equation}\label{eq:spectral-growth-multiplicity}
    \lambda_{h,j}\geq \widetilde{c}j^{2/d},\qquad j=1,\ldots,N.
\end{equation}
This assumption is a discrete counterpart of Weyl's law~\cite{Weyl},
with its corresponding bounds for quasi-uniform meshes discussed in~\cite{XuZhangZikatanov2017}. 
Let the distinct eigenvalues be 
\begin{equation}\label{eq:distinct-spectrum}
  0<\widetilde\lambda_1<\cdots<\widetilde\lambda_{\widetilde N}, \qquad \widetilde{N}\leq N,
\end{equation}
and let $\widetilde P_j$ be the orthogonal projector onto the corresponding eigenspace. Define $b_j=\widetilde P_jb$ for $j=1,\ldots,\widetilde N$.
Since the $j$th distinct eigenvalue is no smaller than the $j$th eigenvalue counted with multiplicity, we have
\begin{equation}\label{eq:spectral-growth}
  \widetilde\lambda_j\ge \widetilde{c}\, j^{2/d},
  \qquad j=1,\dots,\widetilde N.
\end{equation}
We make the following assumption on the observable. 
\begin{assumption}\label{a:observable}
For a fixed integer $1\leq r<\widetilde N$,
 assume that the observable~\eqref{y} satisfies
\begin{equation} 
  b_j=\widetilde{P}_jb\ne 0,\qquad j=1,\ldots,r,
  \label{eq:visible-modes}
\end{equation}
 where $\widetilde{P}_j$ is the orthogonal projection onto the eigenspace corresponding to the distinct eigenvalue $\widetilde{\lambda}_j$.
\end{assumption}
In other words, Assumption~\ref{a:observable} states that the observable $y$ ``sees'' the first $r$ distinct eigenvalues. For a grid-point observation $b=e_i$, the condition $b_j\ne0$ fails precisely when all vectors in the eigenspace corresponding to $\widetilde\lambda_j$ vanish at the observed node. More generally, $b_j=0$ if and only if $b$ is orthogonal to that eigenspace.

Let
$\widetilde E_r$ be the span of the visible spectral components \eqref{eq:visible-modes} and $P_r$
the orthogonal projection onto it:
\begin{equation}
  \widetilde E_r:= \spanop\{b_1,\ldots,b_r\},\qquad
  P_r:= \sum_{j=1}^{r}\frac{b_jb_j^T}{\|b_j\|_2^2}.
\end{equation}
By Assumption~\ref{a:observable}, the $b_j$'s are nonzero and mutually
orthogonal because they lie in eigenspaces corresponding to distinct eigenvalues. Hence,
$\dim \widetilde E_r=r$ and $P_r$ is an orthogonal projection that commutes with
$A_h$, and consequently, with $\Phi_\tau$.  

\begin{lemma}[Semigroup tail decay]\label{lem:tail-decay}
For any \(k\ge 0\),
\begin{equation}
  \left\|(I-P_r)\Phi_\tau^k b\right\|_2
  \le
  e^{-k\tau\widetilde\lambda_{r+1}} \|(I-P_r)b\|_2
  \le
  e^{-\widetilde{c} k\tau (r+1)^{2/d}}\|b\|_2.
  \label{eq:tail-decay}
\end{equation}
\end{lemma}

\begin{proof}
Since $I-P_r$ commutes with $\Phi_\tau$ and the components $b_j$ are
mutually orthogonal,
\[
  \left\|(I-P_r)\Phi_\tau^k b\right\|_2^2
  =
  \sum_{j=r+1}^{\widetilde N} \|b_j\|_2^2\, e^{-2k\tau\widetilde\lambda_j}
   \le
  e^{-2k\tau\widetilde\lambda_{r+1}}
  \sum_{j=r+1}^{\widetilde N} \|b_j\|_2^2.
\]
Taking square roots and using \eqref{eq:spectral-growth} proves the stated inequalities.
\end{proof}

\begin{remark}
 For $k\geq 1$ and a given tolerance $0<\varepsilon<1$, a sufficient condition for the relative tail bound to be at most $\varepsilon$ is
 \begin{equation}
  e^{-\widetilde{c} k\tau (r+1)^{2/d}}\leq \varepsilon,
  \qquad\text{or equivalently,}\qquad
  (r+1)^{2/d}\geq  \frac{\log(1/\varepsilon)}{\widetilde{c} k\tau}.
\end{equation}
Thus a sufficient spectral cutoff is of order
$(\log(1/\varepsilon)/(\widetilde{c}k\tau))^{d/2}$ and does not
explicitly grow with the mesh dimension. In other words, a sufficient
number of delay terms $(r+1)$ to ensure an accuracy requirement
$\varepsilon$ depends explicitly on the spectral decay
\eqref{eq:spectral-growth-multiplicity}  of the semigroup
$\Phi_\tau$, rather than on its dimension $N$. This is consistent with the
numerical evidence of FML.
Note that this is an upper bound obtained by modal truncation, not an asymptotic estimate.
\end{remark}

\subsection{Smoothed Krylov Space}
The preceding lemma shows that the sufficient spectral cutoff decreases as the smoothing interval $k\tau$ increases. Such a smoothing interval, often called burn-in, is frequently used in FML implementations~\cite{Churchill_2023}. We now introduce the \emph{smoothed Krylov space}:
\begin{definition}[Smoothed semigroup Krylov space]
For integers $s\geq0$ and $m\geq1$, define
\begin{equation} 
  \mathcal{K}_m^{(s)}(\Phi_\tau,b) \coloneqq
  \spanop\{\Phi_\tau^s b,\Phi_\tau^{s+1}b,\ldots,\Phi_\tau^{s+m-1}b\},
  \qquad s\ge 0.
  \label{smooth-krylov}
\end{equation}
For $s=0$, it becomes the unshifted observable semigroup Krylov space \eqref{Krylov}.
\end{definition}
%

Let 
\begin{equation} \label{mu}
  \widehat b_j\coloneqq \frac{b_j}{\|b_j\|_2},\qquad
  \beta_j\coloneqq \|b_j\|_2>0,\qquad
  \mu_j\coloneqq e^{-\tau\widetilde\lambda_j},\qquad j=1,\ldots,r,
\end{equation}
so that the vectors $\widehat b_j$ form an orthonormal basis of
$\widetilde E_r$. Since $b_j$ lies in the eigenspace of $\widetilde\lambda_j$, we have
$\Phi_\tau b_j=\mu_j b_j$. Hence, for $s\ge0$ and $k=0,\ldots,r-1$,
\begin{equation}
  P_r\Phi_\tau^{s+k} b
  =
  \sum_{j=1}^r \beta_j \mu_j^{s+k} \widehat b_j.
  \label{eq:coeff-expansion}
\end{equation}
Let $X_s\in\mathbb{R}^{r\times r}$ be the matrix whose rows are the coefficients of
$P_r\Phi_\tau^{s}b,\ldots,P_r\Phi_\tau^{s+r-1}b$ in the basis
$\{\widehat b_1,\ldots,\widehat b_r\}$. Then
\begin{equation}
 X_s=YD_s,\qquad
    Y_{k+1,j}=\mu_j^k,\qquad
    D_s=\diag(\beta_1\mu_1^s,\ldots,\beta_r\mu_r^s).
  \label{eq:X-vandermonde}
\end{equation}
\begin{lemma} \label{lem:projected-spanning}
If Assumption~\ref{a:observable} holds, then, for every $s\ge 0$,
\begin{equation}
  P_r \mathcal{K}_r^{(s)}(\Phi_\tau,b)=
  \spanop\{P_r\Phi_\tau^s b,\ldots,P_r\Phi_\tau^{s+r-1}b\}=\widetilde E_r.
  \label{eq:projected-krylov-spans-Er}
\end{equation}
\end{lemma}

\begin{proof}
The distinct real eigenvalues $\widetilde\lambda_j$ give distinct
positive nodes $\mu_j$ as defined in~\eqref{mu}. Hence $Y$ is an invertible Vandermonde matrix. The diagonal matrix $D_s$ is also invertible because $\beta_j\ne0$ and $\mu_j>0$. Consequently $X_s$ is invertible, and the projected Krylov vectors span $\widetilde E_r$.
\end{proof}
Note that, in general, $\mathcal{K}_r^{(s)}(\Phi_\tau,b)$ itself is
not equal to $\widetilde E_r$ because its vectors contain all the
spectral components.

\begin{theorem}[Approximation by smoothed semigroup Krylov space]\label{thm:smoothed-approximation}
Under Assumption~\ref{a:observable}, for integers $s,k\geq0$, set
\begin{equation}\label{eq:interp-coefficients}
    v_k=(\mu_1^k,\ldots,\mu_r^k)^T,\qquad
    a^{(k)}=Y^{-T}v_k,\qquad
    \Gamma_{r,k}=1+\|a^{(k)}\|_1.
\end{equation}
Then
\begin{equation}\label{eq:distance-to-smoothed-krylov}
\begin{split}
    \dist\bigl(\Phi_\tau^{s+k}b,\mathcal{K}_r^{(s)}(\Phi_\tau,b)\bigr)
    &\leq\Gamma_{r,k}e^{-s\tau\widetilde\lambda_{r+1}}\|(I-P_r)b\|_2\\
    &\leq\Gamma_{r,k}e^{-\widetilde{c}s\tau(r+1)^{2/d}}\|b\|_2.
\end{split}
\end{equation}
Moreover,
\begin{equation}\label{eq:interpolation-factor}
    \Gamma_{r,k}\leq1+\sqrt r\,\|Y^{-T}v_k\|_2
    \leq1+r\|Y^{-1}\|_2.
\end{equation}
In particular, the interpolation coefficients and factor are independent of the smoothing index $s$ and of the nonzero observation amplitudes $\beta_j$.
\end{theorem}
\begin{proof}
By Lemma~\ref{lem:projected-spanning}, there are unique coefficients $a_0^{(k)},\ldots,a_{r-1}^{(k)}$ such that
\begin{equation}\label{eq:low-mode-representation}
    P_r\Phi_\tau^{s+k}b
    =\sum_{\ell=0}^{r-1}a_\ell^{(k)}P_r\Phi_\tau^{s+\ell}b.
\end{equation}
In the basis $\{\widehat b_j\}$, this system is $X_s^Ta^{(k)}=D_sv_k$. Since $X_s^T=D_sY^T$, cancellation of the invertible diagonal factor gives $Y^Ta^{(k)}=v_k$. This proves~\eqref{eq:interp-coefficients} without introducing a smoothing-dependent inverse factor.

Subtracting the corresponding unprojected combination yields
\begin{equation}\label{eq:residual-tail-identity}
    \Phi_\tau^{s+k}b-\sum_{\ell=0}^{r-1}a_\ell^{(k)}\Phi_\tau^{s+\ell}b
    =(I-P_r)\Phi_\tau^{s+k}b
    -\sum_{\ell=0}^{r-1}a_\ell^{(k)}(I-P_r)\Phi_\tau^{s+\ell}b.
\end{equation}
The triangle inequality and Lemma~\ref{lem:tail-decay} give
\begin{equation}\label{eq:IPr}
\begin{split}
    \left\|\Phi_\tau^{s+k}b-\sum_{\ell=0}^{r-1}a_\ell^{(k)}\Phi_\tau^{s+\ell}b\right\|_2
    &\leq\left(1+\sum_{\ell=0}^{r-1}|a_\ell^{(k)}|\right)
    e^{-s\tau\widetilde\lambda_{r+1}}\|(I-P_r)b\|_2.
\end{split}
\end{equation}
The approximation vector belongs to $\mathcal{K}_r^{(s)}$. Finally, $\|a^{(k)}\|_1\leq\sqrt r\|a^{(k)}\|_2$ and $\|v_k\|_2\leq\sqrt r$ because $0<\mu_j<1$. These facts prove~\eqref{eq:distance-to-smoothed-krylov} and~\eqref{eq:interpolation-factor}.
\end{proof}

\begin{remark}[Conditioning and mesh dependence]\label{rem:conditioning}
A direct estimate of $a^{(k)}=X_s^{-T}D_sv_k$ also gives a valid, but weaker, bound
\[
    \Gamma_{r,k}\leq1+\sqrt r\,\|X_s^{-1}\|_2\|b\|_2.
\]
Estimating the factors separately obscures the exact cancellation in the proof. The factor in~\eqref{eq:interpolation-factor} removes this artificial dependence on $s$, but it still depends on the retained semigroup nodes. If these nodes become close, the Vandermonde matrix can be ill-conditioned. 
The theorem concerns exact arithmetic and does not assert robustness of explicitly inverting a Vandermonde matrix.
\end{remark}

\subsection{Closure after Smoothing}\label{sec:burn-memory}

Set $d_\ell=a_\ell^{(r)}$ for $\ell=0,\ldots,r-1$, and reverse the coefficients as in~\eqref{eq:coefficient-reversal}. The following consequence makes the connection to an actual delay model explicit.
\begin{corollary}[Observable closure error bound]\label{cor:smoothed-closure}
Under the conditions of Theorem~\ref{thm:smoothed-approximation}, define $c_\ell=d_{r-1-\ell}$. For every integer $s\geq0$ and every $j\geq s+r-1$,
\begin{equation}\label{eq:smoothed-observable-defect}
    \left|y_{j+1}-\sum_{\ell=0}^{r-1}c_\ell y_{j-\ell}\right|
    \leq\Gamma_{r,r}e^{-\widetilde{c}s\tau(r+1)^{2/d}}\|b\|_2\|U_0\|_2.
\end{equation}
\end{corollary}
\begin{proof}
Apply~\eqref{eq:IPr} with $k=r$ and take the inner product with $U_{j-s-r+1}$. The target index becomes $j+1$ and the retained indices are $j-r+1,\ldots,j$. Contractivity~\eqref{eq:state-contraction} bounds the norm of this state by $\|U_0\|_2$.
\end{proof}

This bound controls the recurrence defect evaluated on the exact observable
history. Error accumulation during autonomous rollout also depends on the
stability of the recurrence.

The oldest retained sample at the earliest admissible forecast is taken at time $s\tau$. Hence
\[
    T_{\rm burn-in}=s\tau,\qquad T_{\rm history}=(r-1)\tau.
\]
Burn-in smoothing suppresses unresolved components in these estimates,
providing a quantitative explanation for its use in FML computations;
see~\cite{Churchill_2023}. It is not necessary in every instance, as the initial state or observable may already have little high-frequency content, or the retained space may contain all visible components. 
\section{Computation of Delay Closures}\label{sec:computation}

We describe two ways to compute a finite-delay model. When the discrete
operator $A_h$ is available, a reduced representation of its semigroup
provides the observation vectors needed to determine the coefficients.
Alternatively, the coefficients can be fitted directly from observable
trajectories. These constructions use the same finite-memory
representation, and
neither requires the full state during subsequent prediction. Numerical
results for both approaches are presented in Section~\ref{sec:examples}.

\subsection{Operator-Based Construction}\label{sec:operator-construction}

In practice, \eqref{Ut} is advanced by a numerical scheme,
\begin{equation}\label{tU}
  \Ut_{k+1}=\tPhi_\tau\Ut_k,\qquad \yt_k=b^T\Ut_k,
\end{equation}
where $\tPhi_\tau:\mathbb R^N\to\mathbb R^N$ is an approximation to the
true semigroup $\Phi_\tau$, constructed from a consistent and stable
numerical time-integration scheme. For a large discrete system, the exact
flow map $\Phi_\tau$ is generally not available explicitly, and $\tPhi_\tau$ provides
a practical approximation to its action. We now discuss how to
construct a low-dimensional space for the available
semigroup $\tPhi_\tau$ required by the delay closure.

In particular, consider the backward Euler propagator
\begin{equation}\label{EB}
    \tPhi_\tau=(I+\tau A_h)^{-1}
       =\eta(A_h+\eta I)^{-1},\qquad \eta=\tau^{-1}.
\end{equation}
Since multiplication by a nonzero scalar does not change a Krylov space,
the corresponding space is generated by the shifted resolvent:
\begin{equation}\label{shift-invert-krylov}
    \mathcal K_m\bigl((A_h+\eta I)^{-1},b\bigr)
    =\operatorname{span}\{b,(A_h+\eta I)^{-1}b,\ldots,
                  (A_h+\eta I)^{-(m-1)}b\},\qquad \eta>0.
\end{equation}
The resolvent has the same eigenvectors as $A_h$, with eigenvalues
$(\lambda_{h,j}+\eta)^{-1}$, and favors the low-frequency components relevant
to the heat equation semigroup.

An orthonormal basis $V_m\in\mathbb R^{N\times m}$ for the space
in~\eqref{shift-invert-krylov} is constructed by shift-invert Lanczos,
using repeated solves with $A_h+\eta I$ and orthogonalization; see
\cite{EricssonRuhe1980,Saad1992,HochbruckLubich1997,EshofHochbruck2006}.
The shifted systems can be solved by preconditioned conjugate gradients
using a multilevel preconditioner, such as BPX~\cite{BramblePasciakXu1990}
or algebraic multigrid~\cite{BrandtMcCormickRuge1982,RugeStuben1987}.

Once the basis is available, define
\begin{equation}\label{eq:reduced-generator}
    T_m=V_m^TA_hV_m,\qquad \xi_b=V_m^Tb,\qquad
    \xi_\ell=e^{-\ell\tau T_m}\xi_b.
\end{equation}
The same basis is used at all required sample times:
\begin{equation}\label{eq:power-projection}
    \Phi_\tau^\ell b=e^{-\ell\tau A_h}b
       \approx V_m\xi_\ell,
       \qquad \ell=0,1,\ldots.
\end{equation}
The exponential is evaluated only for the $m\times m$ matrix $T_m$. Its eigendecomposition
can be reused to generate the reduced sequence $\{\xi_\ell\}$ without
repeated operations on the full state space.
%
For a chosen closure order $r$, determine a forward recurrence by the small
least-squares problem
\begin{equation}\label{eq:reduced-closure-fit}
    \widehat d = \argmin_{d\in\mathbb R^r}
    \left\|\xi_r-\sum_{\ell=0}^{r-1}d_\ell\xi_\ell\right\|_2^2.
\end{equation}
Upon reversing the coefficients,
we obtain the backward-looking delay model
$$
\widehat y_{k+1}=\sum_{\ell=0}^{r-1}\widehat c_\ell y_{k-\ell}, \qquad \widehat c_\ell=\widehat d_{r-1-\ell}.
$$
Thus the reduced basis provides an efficient means of computing the
coefficients of an approximate delay closure.
A small fitting residual in~\eqref{eq:reduced-closure-fit} must be
considered together with the semigroup-action error
in~\eqref{eq:power-projection}.

\subsection{Data-Driven Construction}\label{sec:fml}

The second approach estimates the delay-closure coefficients from
observable histories. It is entirely data driven and does not require
any knowledge of the state equation. This is the common approach by FML, c.f. \cite{Churchill_2023}, which usually utilizes DNN. For the linear system considered in this paper, it can be accomplished by least squares regression.
For a scalar observable and a
candidate closure order $r$, set
\begin{equation}\label{eq:fml-linear-model}
    z_k=(y_k,y_{k-1},\ldots,y_{k-r+1})^T,\qquad
    \widehat y_{k+1}=c^Tz_k
       =\sum_{\ell=0}^{r-1}c_\ell y_{k-\ell}.
\end{equation}
For training samples $\{(z_i,y_i^+)\}_{i=1}^{M}$, where $y_i^+$ is the
observation immediately following the history $z_i$, form $X_r$ with
rows $z_i^T$ and $y_{\rm tar}=(y_1^+,\ldots,y_M^+)^T$. The coefficients
are obtained from
\begin{equation}\label{eq:fml-least-squares-solution}
    \widehat c = \argmin_{c\in\mathbb R^r}
       \|X_rc-y_{\rm tar}\|_2^2.
\end{equation}
The minimum-norm solution is $\widehat c=X_r^\dagger y_{\rm tar}$,
where $X_r^\dagger$ denotes the Moore--Penrose inverse. A rank-revealing factorization can be used
when the histories are nearly linearly dependent.

Training trajectories are generated from independent randomized
initial conditions, represented on the spatial mesh, and a few short
segments, also called ``bursts'', are randomly sampled from each
observable trajectory. Each segment supplies
a history of length $r$ and the next observation for the
least-squares fit~\eqref{eq:fml-least-squares-solution}.
Burst-based training data generation is
standard in FML; see~\cite{Churchill_2023}.

For $q\geq 1$ observables at different sensor locations, the
delay closure for $y_k\in\mathbb R^q$ has the following form with
matrix coefficients:
\begin{equation}\label{eq:fml-vector-model}
    \widehat y_{k+1}=\sum_{\ell=0}^{r-1}C_\ell y_{k-\ell},
 \end{equation}
 where the coefficients are computed by
\begin{equation}
    \widehat\Theta = \argmin_{\Theta\in\mathbb R^{qr\times q}}
    \|X_r\Theta-Y_{\rm tar}\|_F^2, \qquad
    \Theta\triangleq (C_0^T,\dots, C_{r-1}^T)^T.
\end{equation}
The scalar observable
formulation \eqref{eq:fml-linear-model} is recovered
when $q=1$.

\section{Numerical Examples}\label{sec:examples}

The numerical examples address three aspects of the preceding analysis:
accurate approximation of semigroup actions, the effective rank of the
resulting observation sequence under mesh refinement, and autonomous
prediction using coefficients fitted directly from data. We first present
the operator-based tests and then the data-driven experiments. The latter
include constant and spatially varying diffusivities, the
placement of sensors, and applications to point-observable
prediction in two and three dimensions.

All examples use the unit hypercube with homogeneous Dirichlet conditions
on the coordinate-zero faces and homogeneous Neumann conditions on the
coordinate-one faces,
\[
    u|_{x_\ell=0}=0,\qquad
    a(x)\partial_{x_\ell}u|_{x_\ell=1}=0,
    \qquad \ell=1,\ldots,d.
\]

\subsection{Semigroup Approximation and Effective Rank}
\label{sec:examples-operator}

The first test measures the relative error of the approximation
in~\eqref{eq:power-projection} at one sampling interval:
\begin{equation}\label{eq:examples-semigroup-error}
    E_m(\tau)=
    \frac{\|V_m e^{-\tau T_m}\xi_b-e^{-\tau A_h}b\|_2}
         {\|e^{-\tau A_h}b\|_2}.
\end{equation}
Table~\ref{tab:semigroup-convergence} collects the constant-
diffusivity results in one dimension ($d=1$), with $N=1,024$, and
three dimensions ($d=3$), with $N=64^3=262,144$. A basis of dimension
$m\ll N$ gives high accuracy in both cases.
At $\tau=1$, $m=10$ gives errors of $9.3\times10^{-11}$ in one dimension
and $5.0\times10^{-10}$ in three dimensions. The shorter interval
$\tau=0.01$ requires a larger basis: the respective errors are
$8.8\times10^{-12}$ at $m=25$ and $4.8\times10^{-10}$ at $m=30$.
Once the error is small, further increases in $m$ do not produce a
monotone improvement. These tests assess semigroup-action approximation,
not the rollout error of a fitted recurrence.

\begin{table}[tb]
\centering
\caption{Relative semigroup-action error $E_m(\tau)$ for $a\equiv1$.
The state dimensions are $N=1024$ in one dimension and $N=64^3$ in three
dimensions.}
\label{tab:semigroup-convergence}
\small
\setlength{\tabcolsep}{3.8pt}
\begin{tabular}{ccrrrrrr}
\toprule
$d$ & $\tau$ & $m=5$ & $m=10$ & $m=15$ & $m=20$ & $m=25$ & $m=30$\\
\midrule
1 & $1.00$ & $4.6\!\times\!10^{-5}$ & $9.3\!\times\!10^{-11}$ & $5.1\!\times\!10^{-11}$ & $2.2\!\times\!10^{-10}$ & $3.7\!\times\!10^{-10}$ & $4.9\!\times\!10^{-11}$ \\
1 & $0.10$ & $5.0\!\times\!10^{-3}$ & $2.7\!\times\!10^{-8}$ & $4.0\!\times\!10^{-12}$ & $1.0\!\times\!10^{-11}$ & $1.1\!\times\!10^{-11}$ & $1.6\!\times\!10^{-11}$ \\
1 & $0.01$ & $1.1\!\times\!10^{-2}$ & $2.5\!\times\!10^{-4}$ & $1.7\!\times\!10^{-6}$ & $9.4\!\times\!10^{-9}$ & $8.8\!\times\!10^{-12}$ & $2.6\!\times\!10^{-12}$ \\
\midrule
3 & $1.00$ & $1.1\!\times\!10^{-3}$ & $5.0\!\times\!10^{-10}$ & $3.6\!\times\!10^{-13}$ & $4.8\!\times\!10^{-13}$ & $6.1\!\times\!10^{-12}$ & $2.1\!\times\!10^{-12}$ \\
3 & $0.10$ & $2.8\!\times\!10^{-2}$ & $4.0\!\times\!10^{-5}$ & $4.4\!\times\!10^{-8}$ & $3.0\!\times\!10^{-10}$ & $3.6\!\times\!10^{-11}$ & $3.5\!\times\!10^{-11}$ \\
3 & $0.01$ & $4.8\!\times\!10^{-2}$ & $1.0\!\times\!10^{-3}$ & $1.4\!\times\!10^{-5}$ & $5.4\!\times\!10^{-7}$ & $2.0\!\times\!10^{-8}$ & $4.8\!\times\!10^{-10}$ \\
\bottomrule
\end{tabular}
\end{table}

We next examine whether the reduced observation sequence becomes more
complex as the spatial discretization is refined; specifically, whether
the effective numerical rank $r_{\rm eff}$ depends on the state
dimension $N$ at a prescribed tolerance. With $m=50$, generate
$\xi_\ell=e^{-\ell\tau T_m}\xi_b$ for $\ell=0,\ldots,80$ and form
\begin{equation}\label{eq:examples-numerical-rank}
    W=[\xi_0,\ldots,\xi_{80}],\qquad
    r_{\rm eff}=\#\{j:\sigma_j(W)>10^{-9}\sigma_1(W)\}.
\end{equation}
For this study, the shift is $\eta=5/\tau$ and a point sensor is placed
near $(0.75,\ldots,0.75)$. Table~\ref{tab:mesh-rank-selected} reports the
constant-coefficient results. Here $n$ denotes the number of retained
unknowns per coordinate, so $h=1/n$ and $N=n^d$.

\begin{table}[tb]
\centering
\caption{Effective rank under spatial refinement for $a\equiv1$.
The numerical-rank threshold is $10^{-9}$, the reduced basis size is
$m=50$, and the sequence contains $81$ columns. $n$ is the number of
retained unknowns per coordinate, with $N=n^d$.}
\label{tab:mesh-rank-selected}
\small
\begin{tabular}{crrrrr}
\toprule
 & & & \multicolumn{3}{c}{$r_{\rm eff}$}\\
\cmidrule(lr){4-6}
$d$ & $n$ & $N$ & $\tau=0.01$ & $\tau=0.1$ & $\tau=1$\\
\midrule
1 & 256 & $256$ & 9 & 4 & 2 \\
1 & 1024 & $1\,024$ & 9 & 4 & 2 \\
1 & 4096 & $4\,096$ & 8 & 4 & 2 \\
\midrule
2 & 256 & $65\,536$ & 11 & 5 & 2 \\
2 & 512 & $262\,144$ & 11 & 5 & 2 \\
2 & 1024 & $1\,048\,576$ & 10 & 4 & 2 \\
\midrule
3 & 16 & $4\,096$ & 12 & 5 & 2 \\
3 & 32 & $32\,768$ & 11 & 5 & 2 \\
3 & 64 & $262\,144$ & 11 & 4 & 2 \\
3 & 128 & $2\,097\,152$ & 11 & 4 & 2 \\
3 & 256 & $16\,777\,216$ & 10 & 4 & 2 \\
\bottomrule
\end{tabular}
\end{table}

The effective ranks remain small throughout the tested mesh ranges
and show little dependence on $N$. In three dimensions, increasing $N$
from $4,096$ to $16,777,216$ changes the reported rank from $12$ to $10$
at $\tau=0.01$; at $\tau=1$, it remains $2$. The one- and two-dimensional
tests show similarly small variations in effective rank under grid refinement. The larger ranks at smaller $\tau$ are
consistent with weaker damping of high spectral modes over a shorter sampling interval.

\FloatBarrier

\subsection{Data-Driven Delay Closure Predictions}
\label{sec:examples-data-driven}

We now fit the observable delay models \eqref{eq:fml-linear-model} of Section~\ref{sec:fml} directly
from data.
We assess
their autonomous prediction accuracy using full-model solutions as
references and report mean absolute errors (MAEs). We consider constant
diffusivity $a(x)\equiv1$ and spatially varying diffusivity, with the
coefficient in each cell drawn from a uniform distribution on $(1,2)$
and held fixed throughout training and testing.

\begin{figure}[tb]
\centering
\includegraphics[width=0.45\linewidth]
{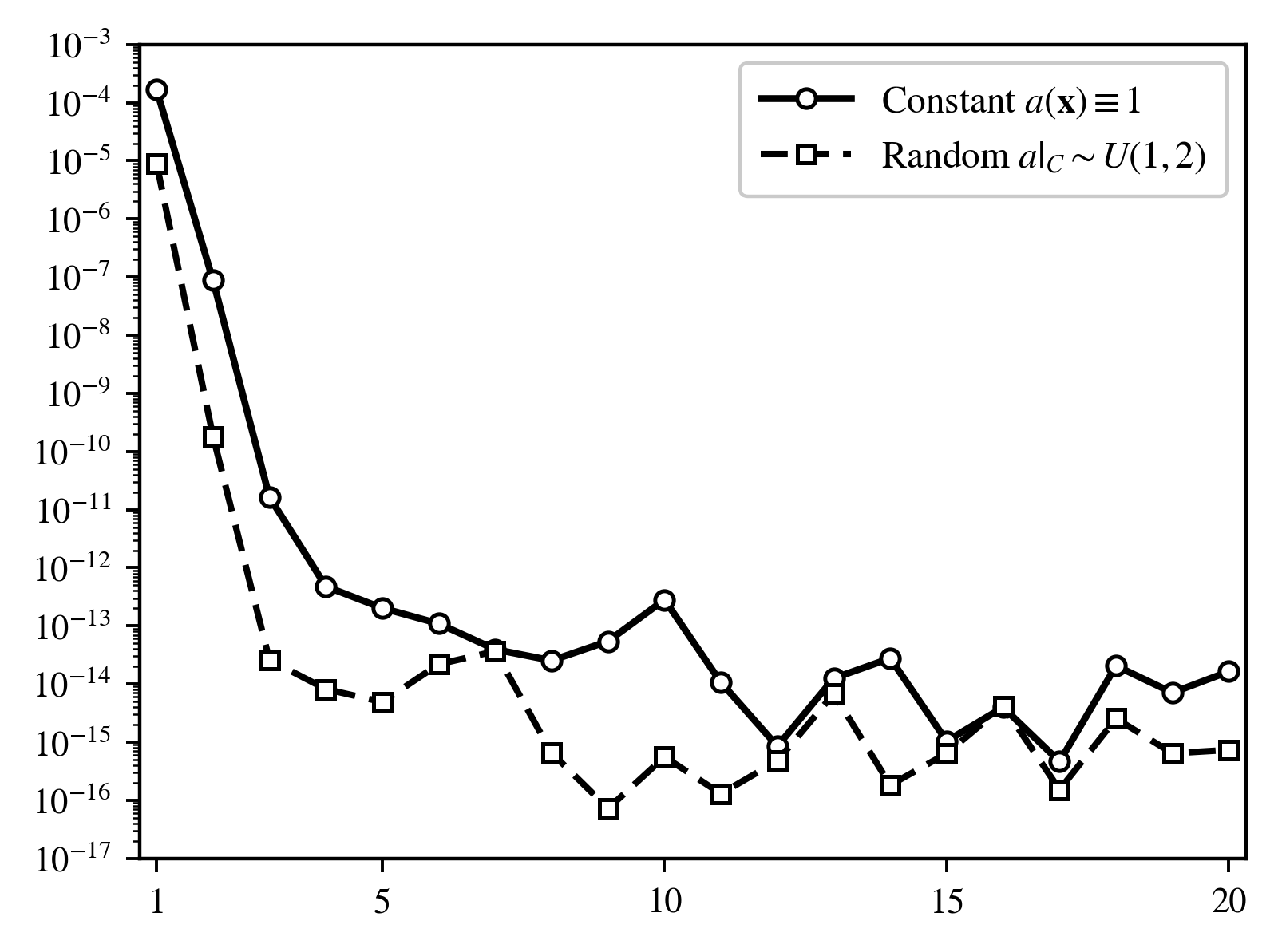}
\includegraphics[width=0.45\linewidth]
{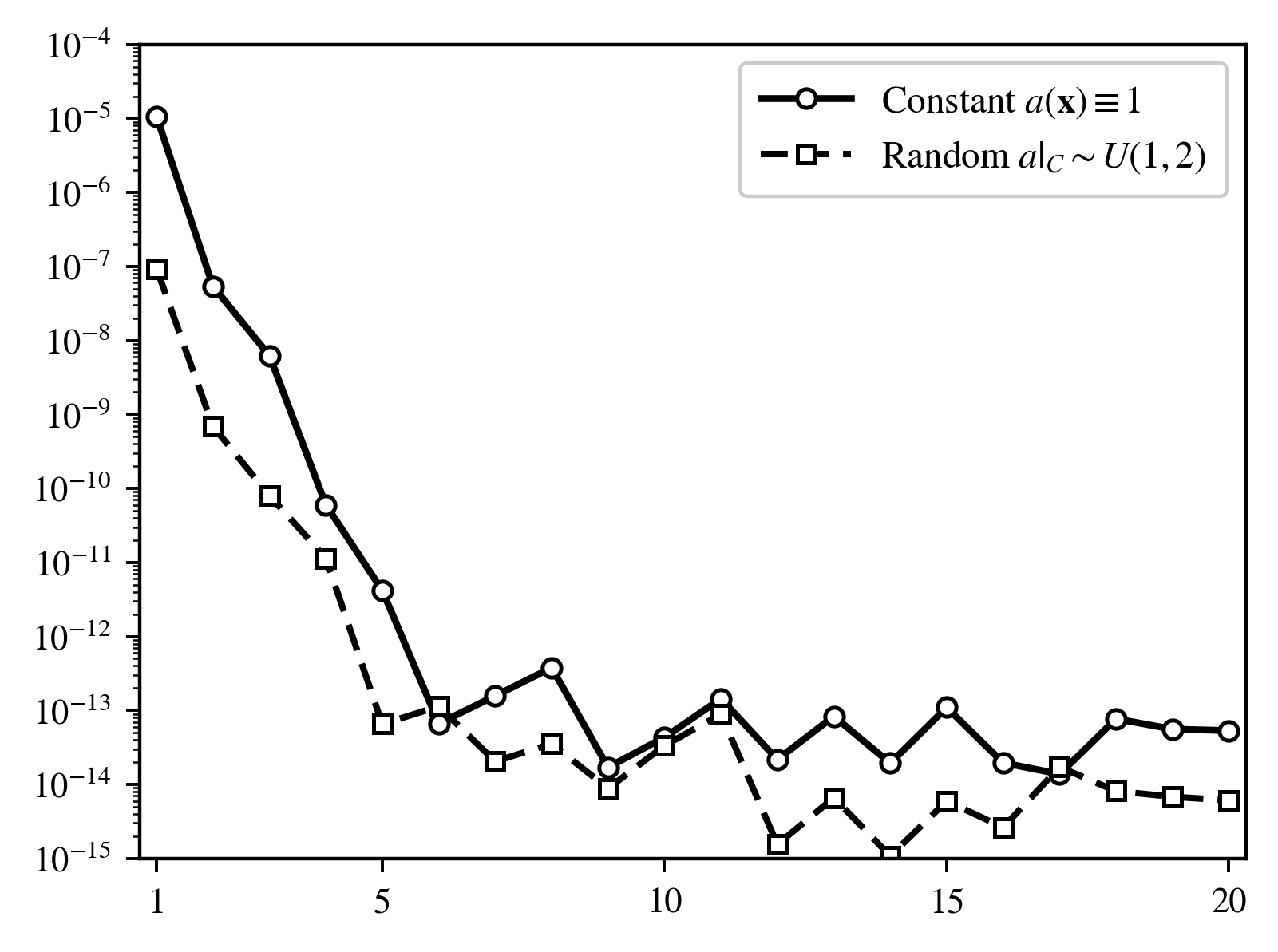}
\includegraphics[width=0.8\linewidth]
{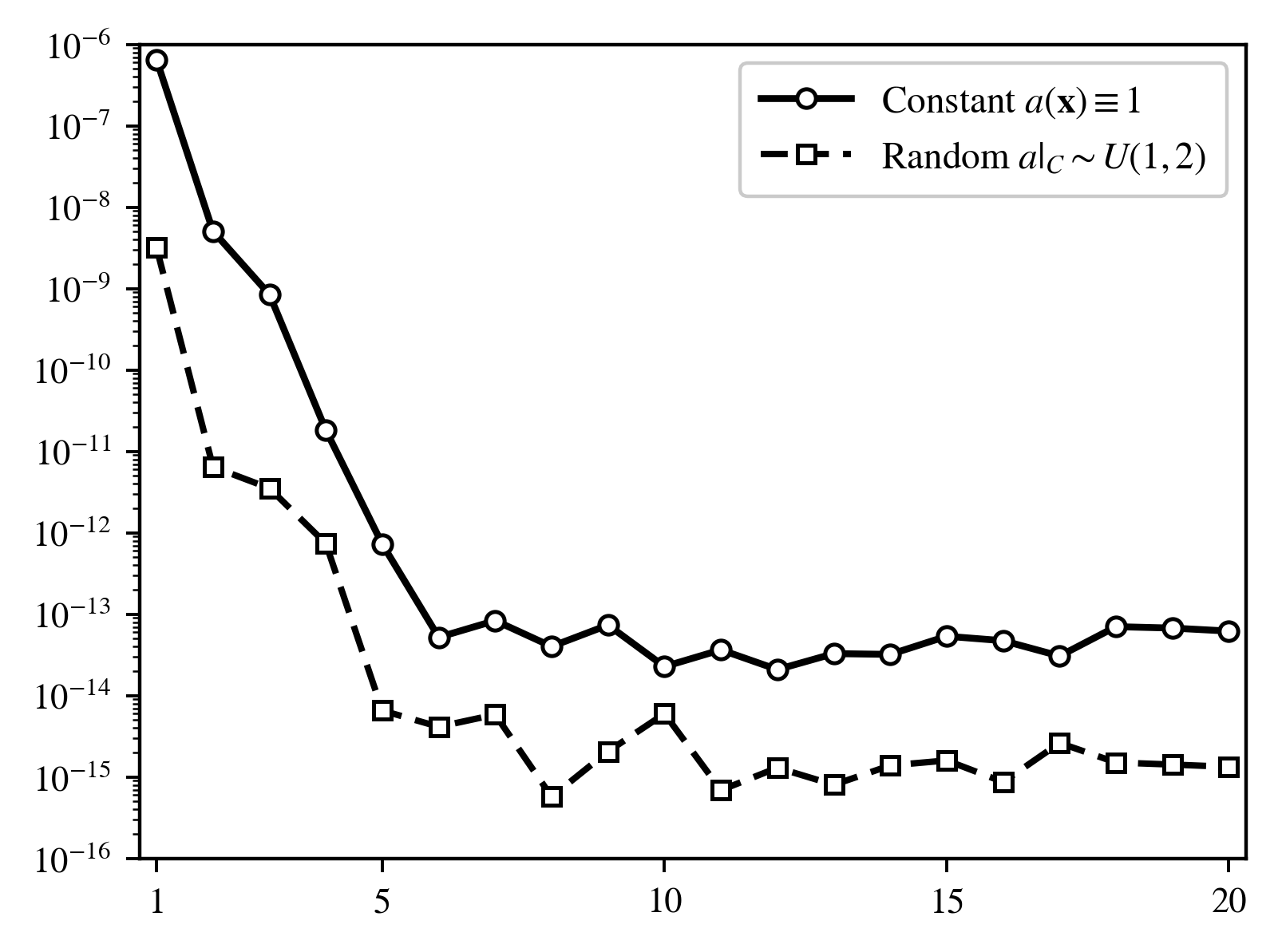}
\caption{Data-driven scalar delay closure in one, two, and three
spatial dimensions (top left, top right, and bottom, respectively).
Reported mean rollout MAE versus delay order $r$ for predictions through
$T=1$ with sampling interval $\tau=0.01$. In each dimension the sensor
is at $x^*=(0.75,\ldots,0.75)$. Solid curves (circles) correspond to
$a\equiv1$; dashed curves (squares) correspond to a fixed realization
of cellwise diffusivities drawn uniformly from $(1,2)$. Each case uses
$20$ independent test trajectories.}
\label{fig:scalar}
\end{figure}

We first consider a single observable in one, two, and
three dimensions, with the sensor at $x^*=(0.75,\ldots,0.75)$. We use
\eqref{eq:fml-linear-model} with time step $\tau=0.01$ to predict the
observable for initial conditions not included in the training data.
Figure~\ref{fig:scalar} shows the reported rollout errors through
$T=1$ as functions of the delay order $r$. In both the constant- and
variable-diffusivity cases, the errors decrease rapidly as $r$ increases,
reaching approximately $10^{-13}$--$10^{-14}$ with about ten delay terms.
These results are consistent with the low-rank semigroup analysis.
The delay model~\eqref{eq:fml-linear-model} uses no solution values at
other spatial locations during prediction; it is a scalar recurrence
whose evaluation cost does not scale with the full state dimension $N$.

\begin{figure}[tb]
\centering
\begin{minipage}[t]{0.47\linewidth}
\centering
\includegraphics[width=\linewidth]{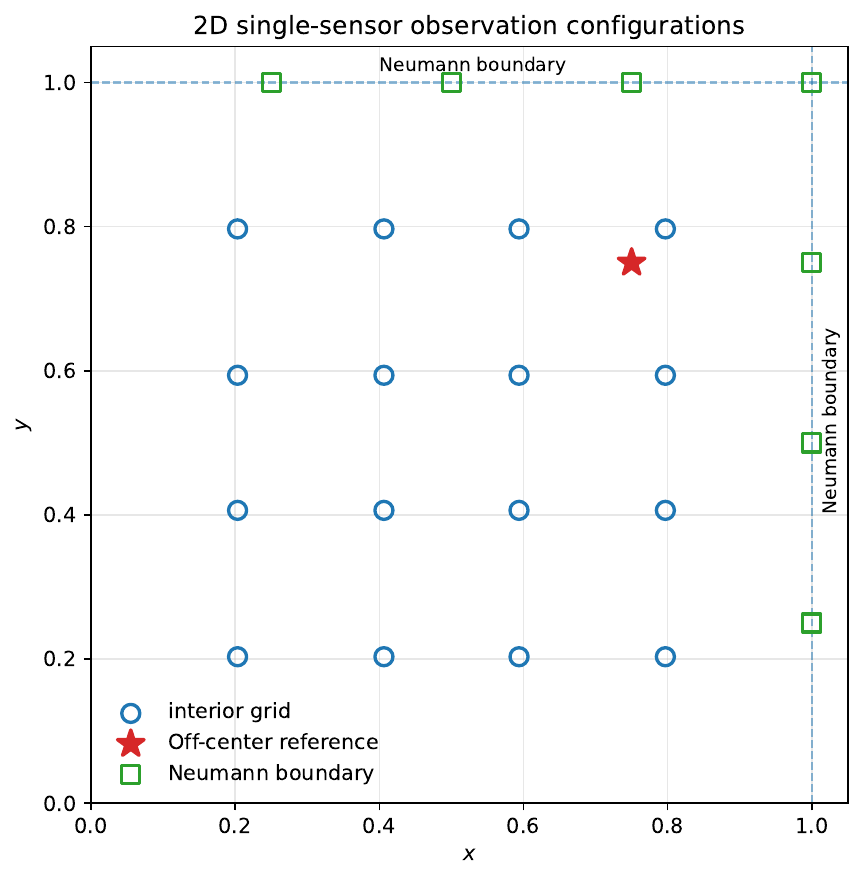}
\par\smallskip
{\small (a) Two dimensions: 24 sensor locations.}
\end{minipage}\hfill
\begin{minipage}[t]{0.49\linewidth}
\centering
\includegraphics[width=\linewidth]{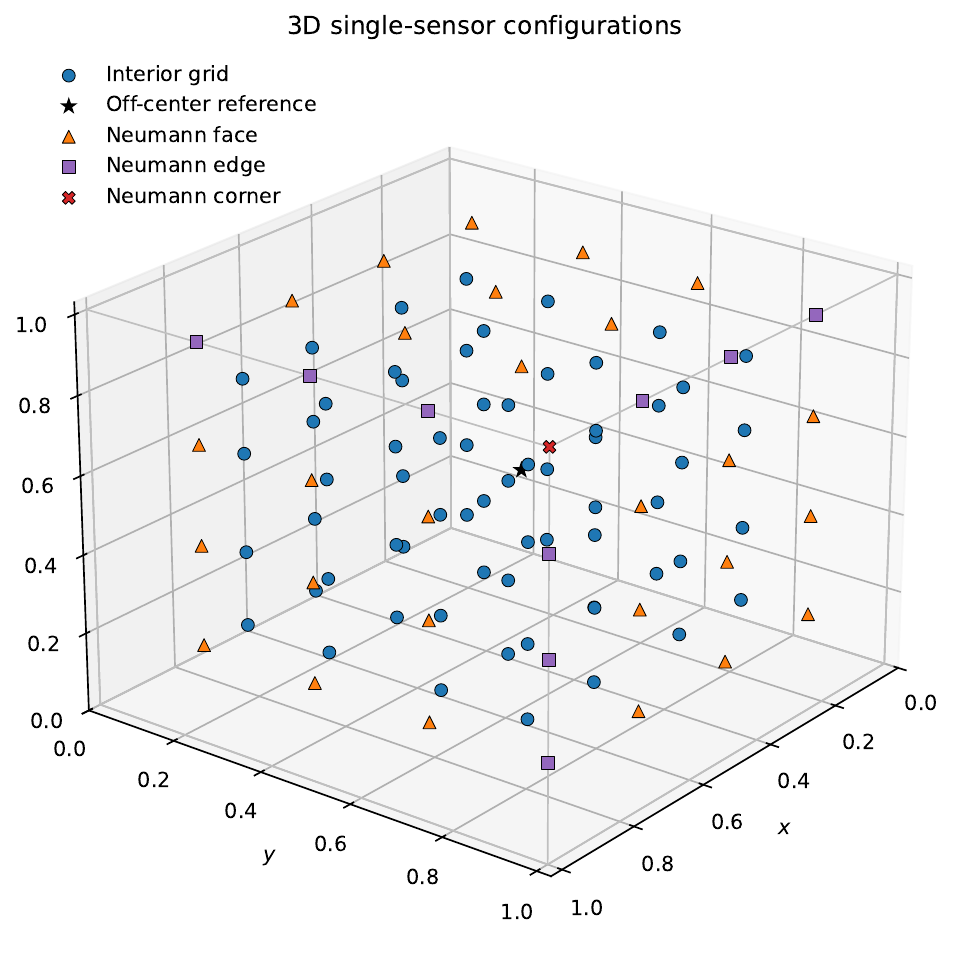}
\par\smallskip
{\small (b) Three dimensions: 102 sensor locations.}
\end{minipage}
\caption{Sensor locations for the single-observable experiments:
(a) $24$ points in two dimensions and (b) $102$ points in three
dimensions. Each point defines a separate scalar-observable model.
The sets include interior locations and points on homogeneous Neumann
boundaries; Dirichlet-boundary locations are excluded. Stars indicate
the reference points $(0.75,0.75)$ and $(0.75,0.75,0.75)$.}
\label{fig:multidimensional-sensor-layouts}
\end{figure}

We next consider constant diffusivity in two ($d=2$) and three
($d=3$) dimensions, varying the location of a single point observable.
Figure~\ref{fig:multidimensional-sensor-layouts} shows the sensor
locations. In two dimensions, the study includes $16$ interior locations
from $\{0.2,0.4,0.6,0.8\}^2$, seven locations on Neumann boundaries,
and the reference location $(0.75,0.75)$ (marked by a star), for a total
of $24$ locations. In three dimensions, the corresponding set comprises
$64$ interior locations from $\{0.2,0.4,0.6,0.8\}^3$, $37$ locations
on Neumann boundaries, and the reference point $(0.75,0.75,0.75)$, for
a total of $102$ locations.

\begin{figure}[htbp]
\centering
\includegraphics[width=0.8\linewidth]{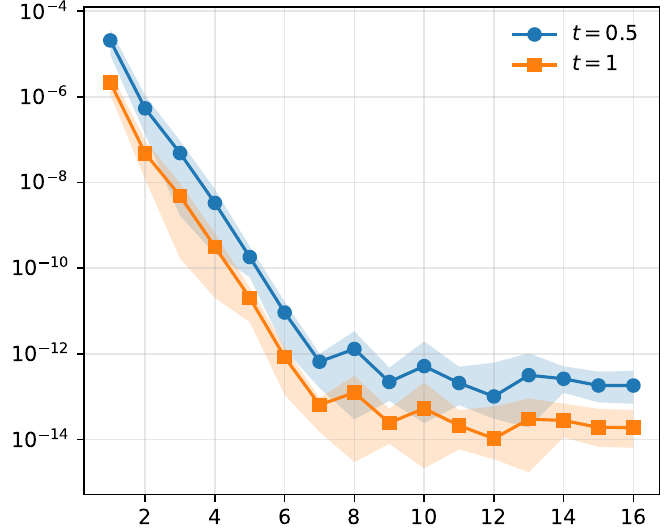}
\par\smallskip
{\small (a) Two dimensions.}
\par\medskip
\includegraphics[width=0.8\linewidth]{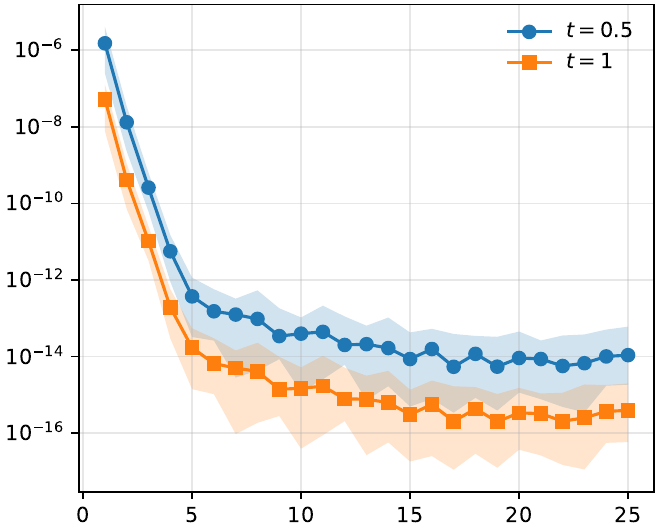}
\par\smallskip
{\small (b) Three dimensions.}
\caption{Autonomous prediction errors for scalar point observables
in (a) two and (b) three dimensions. Curves show terminal-time rollout
MAE at $t=0.5$ and $t=1$ versus the delay order $r$ (denoted by $p$ in
the plots). At each sensor, MAE is computed over $50$ independent test
trajectories; the solid curves average these errors across the sensor
locations shown in Figure~\ref{fig:multidimensional-sensor-layouts}.
Shaded bands indicate the minimum--maximum range across locations, not
confidence intervals. The horizontal axis shows the delay order $r$
on a linear scale, while MAE is displayed on a logarithmic scale.}
\label{fig:multidimensional-rollout}
\end{figure}

Using randomized initial conditions, we construct separate scalar
delay models~\eqref{eq:fml-linear-model} for each sensor location,
with sampling interval $\tau=0.01$. We evaluate autonomous predictions
through $T=1$ against reference solutions from the full state equation.
Figure~\ref{fig:multidimensional-rollout} shows the errors at $t=0.5$
and $t=1$ versus delay order $r$. The solid curves give the mean MAE
across sensor locations, and the shaded regions show the minimum and
maximum. In both dimensions, the error decreases rapidly over the first
several delay orders and then reaches a high-accuracy plateau.
In two dimensions, every location attains an absolute error of
$10^{-12}$ at both test times, and $19$ of the $24$ locations attain
$10^{-14}$ at $t=1$, for some $r\leq16$. In three dimensions, all $102$
locations attain $10^{-14}$ at $t=0.5$ and $10^{-15}$ at $t=1$ for
some $r\leq25$. Thus, accurate rollout prediction is obtained across
a wide range of sensor locations.

\FloatBarrier

\section{Summary}

We have studied finite-delay closure for point observables of spatially
discretized heat equations. Exact sampled-time closure follows from the
semigroup Krylov structure, while parabolic smoothing permits accurate
approximation with substantially fewer delay terms than the full state
dimension. The analysis characterizes the influence of spectral decay
and smoothing on the closure error without requiring reconstruction of
the full state. We have considered two complementary computational
approaches: shift-invert Krylov approximation of the semigroup and
direct least-squares fitting from observable trajectories. Numerical
results in one, two, and three spatial dimensions demonstrate small
effective ranks under mesh refinement and accurate observable
prediction with short delay histories for constant and spatially
varying diffusivities. These findings support finite-delay models as
an efficient means of predicting observable dynamics in diffusion
problems.

\bibliographystyle{siamplain}
\bibliography{references}
\end{document}